\documentclass[a4paper,10pt]{article}

\usepackage{latexsym, amsmath, amsfonts,amssymb, amscd}

\usepackage[all,2cell,ps]{xy}

\usepackage{a4wide}
\usepackage[latin1]{inputenc}
\usepackage[T1]{fontenc}

\usepackage[amsmath,thmmarks,thref]{ntheorem}

\usepackage[latin1]{inputenc}
\usepackage{times}
\usepackage[T1]{fontenc}
\usepackage{amsfonts}
\usepackage{amssymb}
\usepackage{amsmath}
\usepackage{fleqn}

\theoremstyle{change}
\theoremheaderfont{\normalfont\bfseries} 
\theorembodyfont{\itshape}
\theoremseparator{.}
\theoremnumbering{arabic}
\theoremsymbol{}

\newtheorem{thm}{Theorem}[section]
\newtheorem{lemma}[thm]{Lemma}

\newtheorem{prop}[thm]{Proposition}

\theorembodyfont{\upshape}

\newtheorem{defn}[thm]{Definition}
\newtheorem{remark}[thm]{Remark}
\newtheorem{example}[thm]{Example}
\newtheorem{notation}[thm]{Notation}
\newtheorem{assumption}[thm]{Assumption}

\theoremstyle{nonumberplain}
\theoremsymbol{\ensuremath{\square}}
\theoremheaderfont{\scshape}
\newtheorem{proof}{Proof}

\newcommand{\refThm}[1]{Theorem \ref{thm : #1}}
\newcommand{\refLemma}[1]{Lemma \ref{thm : #1}}

\newcommand{\refDef}[1]{Definition \ref{thm : #1}}
 
\newcommand{\refSection}[1]{Section \ref{thm : #1}} 
\newcommand{\refRemark}[1]{Remark \ref{thm : #1}}

\newcommand{\refAssumption}[1]{Assumption \ref{thm : #1}}

\newcommand{\labelThm}[1]{\label{thm : #1}}
\newcommand{\labelAssumption}[1]{\label{thm : #1}}
\newcommand{\labelLemma}[1]{\label{thm : #1}}
\newcommand{\labelProp}[1]{\label{thm : #1}}

\newcommand{\labelDef}[1]{\label{thm : #1}}
\newcommand{\labelSection}[1]{\label{thm : #1}}
\newcommand{\labelRemark}[1]{\label{thm : #1}}

\newcommand{\coh}[3]{\cohH^#1(#2,#3)}
\newcommand{\QZ}{\mathbb{Q}/\mathbb{Z}}
\newtheorem{proofLater}{Proof of \refThm{MainThm}}

\newcommand{\cyclic}[1]{\mathbb{Z}/#1\mathbb{Z}}
\newcommand{\class}[1]{[#1]}
\newcommand{\Gal}[2]{ {\mathop \textnormal{Gal}\nolimits (#1/#2) }}

\newcommand{\Aut}{\mathop \textnormal{Aut}\nolimits}

\newcommand{\finiteField}[1]{\mathbb{F}_{#1}}
\newcommand{\Proj}[1]{{{\mathbb{P}}^1}_{\!\!\!{#1}}}  

\DeclareMathOperator{\Norm}{N}

\newcommand{\Hom}{\mathop \textnormal{Hom}\nolimits}
\newcommand{\degree}{\mathop \textnormal{deg}\nolimits}
\newcommand{\localPowerSeriesField}[2]{{#1}(\!(#2)\!)}
\newcommand{\localPowerSeriesRing}[2]{{#1}[\![#2]\!]}
\newcommand{\spec}[1]{\mathop \textnormal{spec}\nolimits(#1)}

\newcommand{\cohH}{\mathop \textnormal{H}\nolimits}
\newcommand{\dFour}{\mathop {D}\nolimits_4}
\newcommand{\dihedral}[1]{\mathop {D}\nolimits_{#1}}

\title{Liftable $\dFour$\!-Covers}
\author{Louis Hugo Brewis\footnote{\noindent Louis Hugo Brewis \newline Instit\"{u}t f\"{u}r Reine Mathematik \newline  Universit\"{a}t Ulm\newline  Helmholtzstrasse 18\newline    D-89069 Ulm (Germany)}}
\begin{document}
\date{}
\maketitle
\abstract{\noindent Let $k$ be an algebraically closed field of characteristic $p$ and let $G \hookrightarrow \Aut_k(\localPowerSeriesRing{k}{t})$ be a faithful action on a local power series ring over $k$\!. Let $R$ be a discrete valuation ring of characteristic $0$ with residue field $k$\!. One asks, whether it is possible to find a faithful action $G \hookrightarrow \Aut_R(\localPowerSeriesRing{R}{t})$ which reduces to the given action, i.e. a {\sl lift}  to characteristic $0$\!. We show that liftable actions exists in the case that $G = \dFour$ and $p = 2$\!. In fact we introduce a family, the supersimple $\dFour$\!-actions, which can always be lifted to characteristic $0$\!.}

\subsection*{Introduction}
Let $C_k/k$ be a curve over a field $k$ of characteristic $p$ and let $G \hookrightarrow \Aut_k(C_k)$ be a $G$\!-action on $C_k$\!. One says this action {\emph{lifts}} to characteristic $0$ if there exists a local ring $R$ of characteristic $0$ with residue field $k$\!, a smooth $R$\!-curve $C_R/R$ together with a map $G \hookrightarrow \Aut_R(C_R)$ which reduces to the given $G$\!-action on $C_k$\!. This is the {\sl global lifting problem} for the group $G$\!.\\

\noindent Similarly one has the {\sl local lifting problem}: let $G \hookrightarrow \Aut_k(\localPowerSeriesRing{k}{t})$\!. We ask when one can find an embedding $G \hookrightarrow \Aut_R(\localPowerSeriesRing{R}{T})$ reducing to the given one. By considering inertia subgroups one sees that each global lifting problem induces, by localisation and completion at each ramification point, a local lifting problem, and in fact, the \emph{local-global} principle of Green--Matignon \cite{GM} states that these two problems are equivalent. For more information on this see also Bertin--M\'{e}zard \cite{BertinMezard} or Henrio \cite{Henrio}.\\

\noindent It is known that the local lifting problem is very difficult when $p$ divides the order of the group $G$\!. Several results are known. Oort--Sekiguchi--Suwa \cite{OSS} showed that all $G$\!-actions lift, if $G$ is the cyclic group of order $p$\!. Later, Green--Matignon \cite{GM} extended this result to the case $G = \cyclic{p^2}$\!. \\

\noindent For groups which are not cyclic, Bouw--Wewers \cite{BouwWewers} showed that all $G$\!-actions lift if $G$ is the dihedral group of order $2p$\!, where $p$ is an odd prime.  Pagot \cite{Pagot} studied the situation where $G = (\cyclic{p})^2$\! and in particular the case $p = 2$\!. Later Pagot \cite{PagotLiftingPaper} showed that all $G$\!-actions lift, if $G$ is the elementary abelian group of order $4$\!. \\

\noindent A weaker question is whether a group $G$  {admits} {\emph{some}} local action in characteristic $p$ which can be lifted to characteristic $0$\!. Matignon \cite{Matignon} has shown that this is true in the case that $G$ is an elementary abelian $p$\!-group, and later Green \cite{Green} showed a similar result in the case that $G$ is a cyclic $p$\!-group.\\

\noindent In \cite{ProblemsAutoCurves}, Matignon asks what the situation for nonabelian $p$\!-groups is. In this paper we give the first results in this direction by studying the situation in the case that $p = 2$\!. Let $\dFour$ be the dihedral group of order $8$\!. We show  that there do exist examples of local $\dFour$\!-actions on $\localPowerSeriesRing{k}{T}$ which admit lifts to characteristic $0$\!. In fact, we exhibit a family of local $\dFour$\!-actions, the {\sl supersimple} $\dFour$\!-actions, which can always be lifted (\refThm{MainThm}). Furthermore, the local degrees of different of these actions are not bounded from above, i.e. the genera of the respective Katz--Gabber compactifications (see \refRemark{KatzGabberRemark}) of these actions are not bounded from above. This provides some evidence for a conjecture of Chinburg--Guralnick--Harbater, see for instance Chinburg \cite{Chinburg}, which states that all actions of dihedral groups $\dihedral{p^n}$ should lift to characteristic $0$\!.
\subsection*{Overview}
\noindent The first part of this paper is a general overview of the Galois theory that we shall be using. In Section \ref{GalTheoryDFour} we study field extensions with Galois group $\dFour$\!. We also interpret these results in the context of covers of curves. Most notably, we give a method for producing $\dFour$\!-Galois extensions by composing two $\cyclic{2}$\!-Galois extensions and taking the Galois closure. \\

\noindent In \refSection{CohomologicalGaloisTheorySection} we focus on the connection between the Galois theory of cyclic extensions and the theory of group cohomology, following Serre \cite{SerreLocalFields}. We then specialize once again to $\dFour$\!-Galois extensions. The group $\dFour$ has several subgroups of order (and others of index) $2$\!, and therefore, contained in a $\dFour$\!-Galois extension are several $\cyclic{2}$\!-subextensions. We study the Galois theoretic connections between these several subextensions. Finally in \refSection{ArtinSchreierTheorySection} we remind the reader of Artin--Schreier theory and its connection with the cohomological interpretation of cyclic Galois theory. \\

\noindent The Sections \ref{ExampleDFourCover} and  \ref{thm : ClassificationSection} are concerned exclusively with the supersimple $\dFour$\!-Galois extensions. Here we deal exclusively with characteristic $2$\!. As we have already pointed out, it is sufficient to study lifting problems in the local context, and therefore, all fields concerned in these two sections will be local power series fields. First we shall define the notion of a local supersimple $\dFour$\!-extension (\refDef{SupersimpleDefn}), and then we classify them in \refThm{ClassificationCor}. \\

\noindent In the second part of this paper our focus shifts to questions of good reduction. Although our goal is to eventually lift local actions, this part is of a global nature. In \refSection{GoodReductionOfGaloisClosures} we consider a $\dFour$\!-Galois cover $C_3 \rightarrow \Proj{K}$ of smooth projective curves over a $2$\!-adic field $K$\!. Let $C_2$  be the quotient of $C_3$ under a nonnormal subgroup of order $2$\!. Given certain good reduction properties of the intermediate cover $C_2 \rightarrow \Proj{K}$\!, we ask when we can deduce that the curve $C_3$ has potentially good reduction. This culminates in \refThm{SmoothnessThm}.\\

\noindent The assumptions required in \refSection{GoodReductionOfGaloisClosures} are difficult to check, and the purpose of \refSection{MainSection} is to give a method for producing $\dFour$\!-Galois covers of curves which satisfy these assumptions. Finally, we conclude by giving an explicit family of examples. Furthermore, by explicitly studying their reductions to characteristic $2$\!, we see that by localizing and completing these families at their ramification points, we obtain liftable examples of local $\dFour$\!-Galois actions in characteristic $2$\!.\\

\noindent Finally, in \refSection{ProofSection} we prove \refThm{MainThm} which is our main result. \\

{\noindent \textbf{{\ref{thm : MainThm}} Theorem.}
 \sl{All supersimple $\dFour$\!-actions lift to characteristic $0$\!.}\\
}

\noindent Essentially the proof goes as follows. We use \refThm{ClassificationCor} to classify explicitly all local supersimple $\dFour$\!-Galois extensions, and we show that each one is isomorphic to one induced (after localized-completion and reduction) from the examples obtained in \refSection{MainSection}.

\subsection*{Notations}
\noindent Let $k$ be an algebraically closed field of characteristic $2$\!. Let $R_0$ denote the Witt vectors of $k$ and $K_0$ its fraction field. We shall reserve the letter $K$ for a finite field extension of the field $K_0$ and $R$ for the normalisation of $R_0$ inside $K_0$\!. The field $K$ will always be assumed to contain $\sqrt{2}$\!.\\

\noindent If $C$ is a smooth projective curve over a field, then we write $g(C)$ to mean its genus. Lastly, if $A$ is ring, then we write $\Proj{A,z}$ for the projective $A$\!-line with distinguished parameter $z$\!. For the notion of a local degree of different see Serre \cite{SerreLocalFields} Chapter III.\\

\noindent Let $\dFour$ denote the dihedral group of order $8$\!. We fix once and for all two generators $a,b \in \dFour$ with the relations 
\begin{align*}\label{DFour_Generators}
 a^4 = b^2 = 1,\ bab = a^3.                                                                                     \end{align*}

\noindent Whenever $L$ is a field, $G_L$ will denote the absolute Galois group of $L$\!. If $G$ denotes a finite group, then by saying that two $G$\!-Galois extensions of $L$ are isomorphic, we are implicitly assuming that a field isomorphism can be found which respects the identification of the respective Galois groups with $G$\!. \\

\subsection*{Acknowledgements}
\noindent I would like to thank Barry Green, Werner L\"{u}tkebohmert, Stefan Wewers and especially my doctoral supervisor Irene Bouw for their patience, support and valuable input during the period that this work was done and written up. Many improvements have been made over time, and the most important of these are due to their input.
\section{Some Galois theory}
\subsection{General Galois theory of $\dFour$\!-extensions}\label{GalTheoryDFour}
The aim of this section is to state and prove two facts on constructing $\dFour$\!-Galois extensions.  We start by studying the situation for field extensions. Later we shall also interpret the results in the context of covers of smooth algebraic curves. \\

\noindent Let $L_0$ be a field. We assume that we are given two $\cyclic{2}$\!-Galois extensions $ L_0 \subset L_1$ and $ L_1 \subset L_2$\!. The extension $L_0 \subset L_2$ is of degree $4$\!, but not necessarily Galois. We denote the Galois closure of this extension by $L_0 \subset \tilde{L}_2$ and the Galois group by $G$\!. The following lemma will be crucial to our studies later on in this work.
\begin{lemma}\labelLemma{GeneratingDFourExtensions}
  Assume that $L_0 \subset L_2$ is not a Galois extension, i.e. $L_2 \neq \tilde{L}_2$\!. Then the Galois extension $L_0 \subset \tilde{L}_2$ is a $\dFour$\!-extension, i.e. $G \simeq \dFour$\!.
\end{lemma}
\begin{proof}
  \noindent Since $L_0 \subset L_2$ is an  extension of degree $4$\!, we notice that $G \subset S_{4}$\!. We also see that $G$ must be nonabelian since $L_0 \subset L_2$ is not Galois. Furthermore, the fact that $L_0 \subset L_1$ and $L_1 \subset L_2$ are both Galois immediately places restrictions on the subgroups of $G$\!. One checks that all subgroups of $S_4$ satisfying all these conditions are  isomorphic to $\dFour$\!.
\end{proof}
We leave the proof of the following lemma to the reader.
\begin{lemma}\labelLemma{GaloisDFourIdentification}
  Assume the notation of above. Then there exists an isomorphism 
 \begin{align*}
   \Gal{\tilde{L}_2}{L_0} \simeq \dFour  
 \end{align*}
   such that $L_2$ is the fixed field under the subgroup \begin{align*}  \left\langle b \right\rangle \subset \dFour\simeq \Gal{\tilde{L}_2}{L_0}.\end{align*}
\end{lemma}

\begin{remark}
  Notice that the results of this section can also be applied to separable covers of curves. We shall leave the details to the reader.
\end{remark}
\subsection{Cohomological Galois theory of fields}\labelSection{CohomologicalGaloisTheorySection}
In this section we shall gather some more facts on the Galois theory of fields, and in particular its cohomological interpretation. Our reference is essentially the book of Serre \cite{SerreLocalFields}.\\

\noindent Let $L_1/L_0$ be a $\cyclic{2}$\!-Galois extension of the field $L_0$\!. It is known that we have the inflation-restriction exact sequence (see Serre \cite{SerreLocalFields} p.118)
 \begin{align*}
    0 \rightarrow \cohH^1(\Gal{L_1}{L_0},\mathbb{Q}/\mathbb{Z}) \rightarrow  \cohH^1(G_{L_0},\mathbb{Q}/\mathbb{Z}) \rightarrow \cohH^1(G_{L_1},\mathbb{Q}/\mathbb{Z})^{\Gal{L_1}{L_0}} \\
 \rightarrow \cohH^2(\Gal{L_1}{L_0},\mathbb{Q}/\mathbb{Z}) \rightarrow \ldots
 \end{align*}
Furthermore, one sees that 
\begin{align*}
 \cohH^1(\Gal{L_1}{L_0},\mathbb{Q}/\mathbb{Z}) \simeq \cyclic{2}
\end{align*}
and
\begin{align*}
 \cohH^2(\Gal{L_1}{L_0},\mathbb{Q}/\mathbb{Z}) \simeq 0,
\end{align*}
see for instance Serre \cite{SerreLocalFields} p.134. We thus obtain the exact sequence
\begin{align}\label{EasyExactSequence}
    \quad 0 \rightarrow \cohH^1(\Gal{L_1}{L_0},\mathbb{Q}/\mathbb{Z}) \rightarrow \cohH^1(G_{L_0},\mathbb{Q}/\mathbb{Z}) \rightarrow \cohH^1(G_{L_1},\mathbb{Q}/\mathbb{Z})^{\Gal{L_1}{L_0}} \rightarrow 0.
 \end{align}

\noindent Let us very briefly remind ourselves what these cohomology groups mean and how they relate to Galois theory. \\

\noindent Consider the absolute Galois group $G_{L_0}$ of the field ${L_0}$\!. The set of $\cyclic{n}$\!-Galois extensions of $L_0$ corresponds bijectively to the elements of the group $\Hom_{\mathbb{Z}}(G_{L_0},\mathbb{Q}/\mathbb{Z})$\!, and therefore, if we consider the group $\QZ$ as a trivial $G_{L_0}$\!-module, the set of {\emph{cyclic}} Galois extensions of $L_0$ corresponds bijectively to the set of elements of the group  
\begin{align}
  \lim_{n \to \infty}\Hom_{\mathbb{Z}}(G_{L_0},\cyclic{n}) \simeq \Hom_{\mathbb{Z}}(G_{L_0},\mathbb{Q}/\mathbb{Z}) \simeq \cohH^1(G_{L_0},\QZ),
\end{align}
i.e. to the group of $G_{L_0}$\!-{\emph{characters}}.\\

\noindent We can now interpret the exact sequence (\ref{EasyExactSequence}) in terms of Galois theory. Let $L'/L_0$ be a cyclic extension of the field $L_0$\!. This extension corresponds to an element  
  \begin{align*}
    \chi \in \Hom_{\mathbb{Z}}(G_{L_0},\mathbb{Q}/\mathbb{Z}) \simeq \cohH^1(G_{L_0},\mathbb{Q}/\mathbb{Z}).
  \end{align*}
  \noindent Notice that the compositum $L'L_1$ of $L'$ and $L_1$ over $L_0$ is a cyclic extension of $L_1$ of degree dividing $n$\!, and therefore, this corresponds to a character \begin{align*}\chi' \in \Hom_{\mathbb{Z}}(G_{L_1},\mathbb{Q}/\mathbb{Z}) \simeq \cohH^1(G_{L_1},\mathbb{Q}/\mathbb{Z}).\end{align*} One checks that the image of $\chi \in \cohH^1(G_{L_0},\mathbb{Q}/\mathbb{Z})$ under the restriction map 
  \begin{align*}
    \cohH^1(G_{L_0},\mathbb{Q}/\mathbb{Z}) \rightarrow \cohH^1(G_{L_1},\mathbb{Q}/\mathbb{Z})
  \end{align*}
is exactly $\chi'$\!. \\

\noindent Let us now list two properties which will be used later on. We give only a short proof of the last of these and leave the other for the reader.
\begin{lemma}\labelLemma{OrderFourElements}
  Let $\chi \in \cohH^1(G_{L_0},\QZ)$ be an element which maps to an element of order $2$ inside $\cohH^1(G_{L_1},\QZ)$ under the restriction map (\ref{restrictionMap}). Then the order of $\chi$ is a divisor of $4$\!.
\end{lemma}
\begin{proof}
 Use (\ref{EasyExactSequence}).
\end{proof}

\begin{lemma}\labelLemma{ImportantGaloisCor}
  Let $\chi_i \in \cohH^1(G_{L_0},\QZ)$ for $i=1,2$ be two  elements of order $4$ which map to elements of order $2$ inside $\cohH^1(G_{L_1},\QZ)$\!. Then the difference $\chi_1 - \chi_2$ is an element of order at most $2$ inside $\cohH^1(G_{L_0},\QZ)$\!.
\end{lemma}
\begin{proof}
 By assumption and exactness of (\ref{EasyExactSequence}), we see that both $2\chi_1$ and $2\chi_2$ are of order $2$ and in fact contained in the group \begin{align}\cyclic{2} \simeq \Gal{L_1}{L_0} \subset \cohH^1(G_{L_0},\QZ).\end{align} Therefore $2\chi_1 = 2\chi_2$ and the result follows.
\end{proof}

\noindent Let us now apply this formalism of characters to study $\dFour$\!-Galois extensions of a field $L_0$\!. Let $L_0 \subset L$ be $\dFour$\!-Galois and fix an isomorphism $\Gal{L}{L_0} \simeq \dFour$\!. Let $L_1$ be the field $L^{\left\langle a^2,b \right\rangle}$ fixed by the subgroup $\left\langle a^2,b \right\rangle \subset \dFour$\!. Notice that $L_1 \subset L$ is a $(\cyclic{2})^2$\!-Galois extension.\\

\noindent There are exactly three proper subfields of $L$ containing $L_1$ other than $L_1$ itself. These are $ L^{\left\langle a^2 \right\rangle}$\!, $ L^{\left\langle b \right\rangle}$ and $ L^{\left\langle a^2b \right\rangle}$\!. Each is a $\cyclic{2}$\!-Galois extension of $L_1$ and therefore, these fields correspond to order $2$ characters $\chi_{a^2}$\!, $\chi_{b}$ and $\chi_{a^2b}$\!, respectively, of the group $G_{L_1}$\!. Hence we may regard the $\chi_{*}$ as elements of the group
\begin{align}
 \Hom_{\mathbb{Z}}(G_{L_1},\mathbb{Q}/\mathbb{Z}) \simeq H^1(G_{L_1},\mathbb{Q}/\mathbb{Z}).
 \end{align}

\noindent We leave the proof of the following lemma to the reader.
\begin{lemma}\labelLemma{DFourIntermediateLemma}
  We have the following relations.
  \begin{enumerate}
   \item The character $\chi_{a^2}$ is fixed under the Galois action $\Gal{L_1}{L_0}$ on the group $\coh{1}{G_{L_1}}{\QZ}$\!.
   \item The characters $\chi_{b}$ and $\chi_{a^2b}$ are conjugate under this action.
   \item The sum of $\chi_{b}$ and $\chi_{a^2b}$ is $\chi_{a^2}$\!.
  \end{enumerate}
\end{lemma}
 \noindent The following lemma will be useful for lifting $\dFour$\!-actions later on.
 \begin{lemma}\labelLemma{UniqueRemark}
   Let $L_0 \subset L$ and $L_0 \subset L'$ be two $\dFour$\!-Galois extensions. Assume that there exists a $L_0$\!-isomorphism between $L^{\left\langle b \right\rangle}$ and $L'^{\left\langle b \right\rangle}$\!. Then there is also a $L_0$\!-isomorphism between $L$ and $L'$\!.
 \end{lemma}
\begin{proof}
  Use the uniqueness of the Galois closure.
\end{proof}
\begin{remark}
 Although a simple lemma, the above tells us that the essential information of the $\dFour$\!-Galois extension $L_0 \subset L$ is stored inside the subextension $L_0 \subset L^{\left\langle b \right\rangle}$\!.
\end{remark}
\begin{notation}\label{NormNotation}
\noindent From now on, whenever we are given extensions as above and a character \begin{align*}\chi \in \cohH^1(G_{L_1},\QZ),\end{align*} we shall denote by $\Norm \chi$ the norm (some reference refer to this as the {\emph {trace}}) of $\chi$ under the action of $\Gal{L_1}{L_0}$ on $\cohH^1(G_{L_1},\QZ)$\!, i.e. the sum of $\chi$ and its conjugate $\sigma^{*}\chi$\!, where $\sigma$ is the generator of $\Gal{L_1}{L_0}$\!. Furthermore, we reserve the notation $\chi_{a^2}$\!, $\chi_{b}$ and $\chi_{a^2b}$ for the characters corresponding to the $\cyclic{2}$\!-extensions $L^{\left\langle a^2 \right\rangle}$\!, $L^{\left\langle b \right\rangle}$ and $L^{\left\langle a^2b \right\rangle}$ of the field $L_1 = L^{\left\langle a^2,b \right\rangle}$\!. 
\end{notation}
\subsection{Artin--Schreier theory of power series fields in characteristic $2$}\labelSection{ArtinSchreierTheorySection}
\noindent Let $L_z$ be the local field $\localPowerSeriesField{k}{z}$ with parameter $z$\!. The following identification will be used often:
 \begin{align}\label{IdentificationArtinSchreierCohomology}
   \quad\localPowerSeriesField{k}{z}/\wp\localPowerSeriesField{k}{z} \simeq \cohH^1_{et}(\spec{\localPowerSeriesField{k}{z}},\cyclic{2})\simeq \cohH^1(G_{L_z},\QZ)[2]
\end{align}
 where \begin{align*}\wp : \localPowerSeriesField{k}{z} \rightarrow \localPowerSeriesField{k}{z}, \quad y \mapsto y^2 - y\end{align*} is the Artin--Schreier operator in characteristic $2$\!. The identification (\ref{IdentificationArtinSchreierCohomology}) associates to the element $$f \in \localPowerSeriesField{k}{z}$$ the $\cyclic{2}$\!-extension of $\localPowerSeriesField{k}{z}$ generated by $w$\!, where $w$ satisfies 
\begin{align*}
 w^2 - w = f.
\end{align*}
\noindent Often we shall denote the associated class of $f$ simply by \begin{align*}\class{f} \in \cohH^1(G_{L_z},\cyclic{2}) \simeq \cohH^1(G_{L_z},\QZ)[2].\end{align*} Notice that for $f_1$ and $f_2$ both elements of $\localPowerSeriesField{k}{z}$\!, we have that 
\begin{align*}
 \class{f_1} = \class{f_2}
\end{align*}
if and only if there exists a $q \in \localPowerSeriesField{k}{z}$\!, such that 
\begin{align*}
 f_1 = q^2 - q + f_2
\end{align*}
inside the field $\localPowerSeriesField{k}{z}$\!.\\

\noindent  Let $f \in k[z] \subset \localPowerSeriesField{k}{z}$\!. In this case, one can always find a $f_0 \in \localPowerSeriesField{k}{z}$ such that 
  \begin{align*}
    f_0^2 - f_0 = f,
  \end{align*}
  and therefore,
  \begin{align*}
    \class{f} = 0
  \end{align*}
  in this case.  Therefore, if $f:=\sum\limits_{-N \leq i} c_iz^{i} \in \localPowerSeriesField{k}{z}$ is a general element of the field $\localPowerSeriesField{k}{z}$ for $c_i \in k$\!, then 
  \begin{align}\label{ElementExpansion}
    \class{f} = \class{\sum\limits_{-N \leq i} c_iz^{i}} = \class{\sum\limits_{ N \leq i < 0} c_iz^{i}}.
  \end{align}
  Furthermore, we also have
  \begin{align*}
    \class{c_{2m}z^{-2m}} = \class{\sqrt{c_{2m}}z^{-m}}
  \end{align*}
  since $k$ is assumed to be algebraically closed. Therefore, we can also get rid of the terms of $f$ in the expansion (\ref{ElementExpansion}) of degree $\!-2m$\!, where $m$ ranges over the natural numbers.

\subsection{Supersimple $\dFour$\!-extensions}\labelSection{SupersimpleSection}\label{ExampleDFourCover}
We now define and study the the type of $\dFour$\!-extensions that we are interested in lifting. Assume throughout this section that $L_0$ is a local power series field with characteristic $2$\!. 
\begin{defn}\labelDef{SupersimpleDefn}
 A local $\dFour$\!-Galois extension $L_0\subset L$ is said to be {\em supersimple} if the following conditions hold.
 \begin{enumerate}
  \item  The local degree of different of $L^{\left\langle a^2,b \right\rangle} \subset L^{\left\langle a^2\right\rangle}$ is $2$\!,
  \item  the local degree of different of $L_0 \subset L^{\left\langle a^2,b \right\rangle}$ is $2$\!.
 \end{enumerate}
\end{defn}
\begin{remark}\labelRemark{KatzGabberRemark}
   \noindent Let $G$ be a finite $p$\!-group and consider a $G$\!-Galois extension of local power series fields \begin{align*} \localPowerSeriesField{k}{z} \subset L_G.\end{align*} We use $L_G$ with the subscript $G$ to emphasize that we are not restricted to supersimple extensions in this remark. \\

  \noindent It is known that there exists a $G$\!-Galois cover of smooth curves 
  \begin{align*}
    C\rightarrow \Proj{k,z}
  \end{align*}
  which is \'{e}tale over $\mathbb{A}^1_k \subset \Proj{k}$\!, completely branched over the complement  $(z = 0) \in \Proj{k}$\!, and which induces $\localPowerSeriesField{k}{z} \subset L_G$ after localization and completion at $z = 0$\!. This cover is known as the Katz--Gabber cover associated to the extension $\localPowerSeriesField{k}{z} \subset L_G$\!. For details on this and for the more general Katz--Gabber compactification, see for instance the account in Gille \cite{KatzGabberReferences}.\\

  \noindent Applying this to the case $G = \dFour$ with $L_0 = \localPowerSeriesField{k}{z}$ and $L_G = L$\!, one sees that $L_0 \subset L$ is supersimple if and only if \begin{align*}C/\left\langle a^2 \right\rangle \simeq \Proj{k},\end{align*} where $C \rightarrow \Proj{k}$ is the Katz--Gabber cover associated to $L_0 \subset L$\!. Notice that this compactification is therefore a hyperelliptic curve. 
\end{remark}
Let us now construct some examples of supersimple extensions. First we set some notation.
\begin{notation}\label{FieldNotation}
 \noindent From now on, we shall reserve the notation $L_0$ for the local power series field $\localPowerSeriesField{k}{t}$\!, and the notation $L_1$ for the local power series field $\localPowerSeriesField{k}{v}$\!, where the variables $t$ and $v$ are related by \begin{align}\label{tvRelationOne}v^{-2} - v^{-1} = t^{-1}\end{align} Also, we shall let $\sigma$ denote the generator of $\Gal{L_1}{L_0} \simeq \cyclic{2}$\!.
\end{notation} 
\begin{example}
  \noindent In view of \refLemma{GeneratingDFourExtensions}, we now construct some $\cyclic{2}$\!-extensions of $L_1$ which, when considered as degree $4$\!-extensons of $L_0$\!, are not Galois. \\

  \noindent Let $\eta \in k$ and consider the element $f_\eta \in L_1$  given by 
  \begin{align*}
    f_\eta = \eta^2v^{-3} - \eta^2v^{-2} = \eta^2t^{-1}v^{-1}.
  \end{align*}

  \noindent Notice that the sum of $f_\eta$ and its conjugate $\sigma^{*}f_\eta$ is simply $\eta^2t^{-1}$\!. One checks quickly that for $\eta \notin \finiteField{2}$\!, the Artin--Schreier class of $\eta^2t^{-1}$ is non-trivial in the group $\cohH^1(G_{L_1},\QZ)$\!. Therefore, the extension $L_2$ of $L_1$\!, defined by 
  \begin{align*}
    w^2 - w = \eta^2t^{-1}v^{-1},
  \end{align*}
  induces an extension of $L_0$ which is of degree $4$ and not Galois. As we have already pointed out in \refLemma{GeneratingDFourExtensions}, this then produces a $\dFour$\!-Galois extension $L_0 \subset L$ by taking the Galois closure. One notes that, by {\refLemma{GaloisDFourIdentification}} the Galois group can be identified with $\dFour$ such that $L_1$ is the fixed field of $\left\langle a^2,b \right\rangle$ and also $L_2$ that of $\left\langle b \right\rangle$\!. The extension of $L_1$ defined by 
 \begin{align}\label{middleASExtension}
   s^2 - s = \eta^2t^{-1}
 \end{align}
 is, by \refLemma{DFourIntermediateLemma} (3), exactly the field extension defined by $L_1 =  L^{\left\langle a^2,b \right\rangle} \subset  L^{\left\langle a^2 \right\rangle}$\!, which one checks has local degree of different exactly $2$\!. Therefore, the $\dFour$\!-Galois extension $L_0 \subset L$ is supersimple.
\end{example}

\noindent The idea of this example is that it is somewhat representative of  supersimple $\dFour$\!-actions. In fact, it will be useful for classifying them (\refThm{ClassificationCor}). 
\begin{defn}\labelDef{FundDefn}
 For a $\eta \in k$\!, we denote by $\psi_\eta$ the character of $G_{\localPowerSeriesField{k}{v}}$ corresponding to the $\cyclic{2}$\!-extension generated by $w^2 - w = \eta^2t^{-1}v^{-1}$\!, i.e. $\psi_\eta$ denotes the the image of the polynomial $\eta^2t^{-1}v^{-1}$ in $\cohH^1(G_{L_1},\QZ)$ under (\ref{IdentificationArtinSchreierCohomology}). 
\end{defn}
\begin{remark}
 It is important to note that \refDef{FundDefn} depends on the choices of the parameters $t$ and $v$\!.
\end{remark}
\begin{remark}\label{cyclicFourCharacters}
 For the value $\eta = 1$\!, the character $\psi_1$ induces a $\cyclic{4}$\!-Galois extension of $L_0$\!. Furthermore, one can show that there exists a character $\psi_1'$ of the group $G_{L_0}$
 \begin{align*}
   \psi_1' : G_{L_0} \twoheadrightarrow \cyclic{4},
 \end{align*}
 which maps to $\psi_1$ under the restriction mapping
 \begin{align*}
   \Hom_\mathbb{Z}(G_{L_0},\QZ) \simeq \cohH^1(G_{L_0},\QZ) \rightarrow \cohH^1(G_{L_1},\QZ) \simeq \Hom_\mathbb{Z}(G_{L_1},\QZ)
 \end{align*}
 of (\ref{EasyExactSequence}). Furthermore, the character $\psi_1'$ generates the torsion subgroup of 
 \begin{align}
  \cohH^1(G_{L_1},\cyclic{4}) \subset \cohH^1(G_{L_1},\QZ) \simeq \Hom_\mathbb{Z}(G_{L_1},\QZ)
 \end{align}
 of order-\!$4$ characters.
\end{remark}
\begin{remark}\label{SumIdentityRemark}
 Notice that we have the following identity in $\cohH^1(G_{L_1},\QZ)$ for all $\eta \in k$\!.
 \begin{align}\label{SumIdentity}
  \psi_{\eta + 1} = \psi_\eta + \psi_1.
 \end{align}
\end{remark}

\subsection{Classifying supersimple $\dFour$\!-extensions}\labelSection{ClassificationSection}
The aim of this section is to classify the local supersimple $\dFour$\!-Galois extensions. Assume throughout that $\localPowerSeriesField{k}{t} = L_0 \subset L$ is a supersimple $\dFour$\!-Galois extension. 
\begin{lemma}
 By possibly changing the parameter $t$ of $L_0$\!, we may assume that the intermediate field extension $L_0 \subset L^{\left\langle a^2,b \right\rangle}$ is generated by $v$\!, where $v$ and $t$ are related by 
 \begin{align}\label{tvRelationTwo}
   v^{-2} - v^{-1} = t^{-1}.
 \end{align}
\end{lemma}
\begin{proof}
 This result follows from the fact that the local degree of different of $L_0  \subset L^{\left\langle a^2,b \right\rangle}$ is $2$\!.
\end{proof}
From now on we set $L_1 = L^{\left\langle a^2,b \right\rangle}$\!. We consider the elements $t$ and $v$ fixed, and use the notation of \refDef{FundDefn}.
\begin{lemma}\labelLemma{ClassificationLemma}
 \noindent Consider the $\cyclic{2}$\!-Galois extension $L^{\left\langle b \right\rangle} / L_1$ and the associated $G_{L_1}$\!-character \begin{align*}\chi_b \in \cohH^1(G_{L_1},\mathbb{Q},\mathbb{Z})\end{align*} of order $2$\!. Then there exists an $\eta \in k$ such that the $G_{L_1}$\!\!-character $\chi_{b} - \psi_\eta$ is the image of a $2$\!\!-torsion element of $\cohH^1(G_{L_0},\QZ)$ under the restriction map \begin{align*}\cohH(G_{L_0},\QZ) \twoheadrightarrow \cohH(G_{L_1},\QZ).\end{align*}
\end{lemma}
\begin{proof}
 \noindent By definition the local degree of different of $L_1 \subset L^{\left\langle a^2\right\rangle}$ is $2$ and therefore, this extension is generated by an Artin--Schreier equation of the form
  \begin{align}\label{centerExtensionEquation}
    s^2 - s = \alpha v^{-1},
  \end{align}
  for some $\alpha \in k$\!. By \refLemma{DFourIntermediateLemma} the norm of $\chi_b$ is the character $\chi_{a^2}$\!. The latter corresponds to the field extension $L_1 \subset L^{a^2}$ and  therefore corresponds to the Artin--Schreier class $\class{\alpha v^{-1}}$\!.\\

  \noindent Choose $\eta \in k$ such that $\eta^2 + \eta = \alpha$\!. One checks that the Artin--Schreier classes of $\class{\alpha v^{-1}}$ and $\class{\eta^2 t^{-1}}$ are the same inside $\cohH^1(G_1,\QZ)$\!.\\

  \noindent Consider the norm $\Norm \psi_\eta$ of the character $\psi_\eta \in \cohH^1(G_{L_1},\QZ)$\!. We see that this corresponds to the $\cyclic{2}$\!-Galois extension of $L_1$ generated by $\tilde{s}$\!, where $\tilde{s}$ satisfies \begin{align*}\tilde{s}^2 - \tilde{s} = \eta^2 t^{-1}.\end{align*} However, by definition of $\alpha$ and $\eta$\!, this is exactly the extension $L_1 \subset L^{a^2}$\!, see (\ref{centerExtensionEquation}). \\

\noindent Hence the norm $\Norm \psi_\eta$ and the character $\chi_{a^2} = \Norm \chi_{b}$ are equal inside the group $\cohH^1(G_{L_1},\QZ)$\!, and hence the difference $ \chi_{b} - \psi_\eta $ is fixed under the action of $\Gal{L_1}{L_0}$\!.\\

\noindent Therefore, $\chi_b - \psi_\eta$ is an element of $\cohH^1(G_{L_1},\QZ)^{\Gal{L_1}{L_0}}$\!, and thus, by the right exactness of (\ref{EasyExactSequence}), the image of some $\chi' \in \cohH^1(G_{L_0},\QZ)$ under the restriction map \begin{align}\label{restrictionMap} \cohH(G_{L_0},\QZ) \twoheadrightarrow \cohH(G_{L_1},\QZ).\end{align}
\noindent By \refLemma{OrderFourElements}, we may conclude that $\chi'$ has order a divisor of $4$\!.\\

  \noindent By Remark \ref{cyclicFourCharacters}, we notice that $\psi_1$ is also the image of an order $4$ element $\psi_1'$ of $\cohH^1(G_{L_0},\QZ)$\!. Hence, by \refLemma{ImportantGaloisCor}, either $\chi'$ or $\chi' - \psi_1'$ is an order $2$ element of $\cohH^1(G_{L_0},\QZ)$\!. \\

  \noindent If $\chi'$ is of order $2$\!, then we have found a suitable $\eta$ satisfying the hypothesis of the lemma. \\

  \noindent Assume this is not the case, i.e. $\chi' - \psi_1'$ is of order $2$\!. Then the image of $\chi' - \psi_1'$ inside $\cohH^1(G_1,\QZ)$ under the restriction map (\ref{restrictionMap}) is exactly \begin{align*}\chi_b - \psi_\eta - \psi_1 = \chi_b - \psi_{\eta + 1},\end{align*} and therefore, the value $\eta + 1$ satisfies the hypothesis of the lemma.
\end{proof}
\begin{thm}\labelLemma{ClassificationCor}
 There exists a polynomial \begin{align*}Q(t^{-1}) \in k[t^{-1}] \subset \localPowerSeriesField{k}{t}\end{align*} and an $\eta \in k$\!, such that the field extension 
 \begin{align*}L_1 = L^{\left\langle a^2,b \right\rangle} \subset L^{\left\langle b \right\rangle}\end{align*} is generated by an Artin--Schreier equation of the form
 \begin{align}\label{ClassificationASEquation}
   w^2 - w = \eta^2t^{-1}v^{-1} + Q(t^{-1}).
 \end{align}
 Furthermore, the polynomial $Q$ can be chosen to have only odd degree terms in the variable $t^{-1}$\!.
\end{thm}
\begin{proof}
  We let $\eta$ be as in \refLemma{ClassificationLemma}. Let $\chi'$ be an element of $\cohH^1(G_{L_0},\QZ)$ which maps to 
  \begin{align*}
    \chi_b - \psi_\eta \in \cohH^1(G_{L_1},\QZ)
  \end{align*}
  under the restriction map \begin{align*}\cohH(G_{L_0},\QZ) \twoheadrightarrow \cohH(G_{L_1},\QZ)\end{align*} and which has order at most $2$\!.\\

  \noindent The character $\chi'$ corresponds to a cyclic Galois extension of the field $L_0$ of degree at most $2$\!. Therefore, we can find an element $Q$ of the field $L_0 = \localPowerSeriesField{k}{t}$\! with associated Artin--Schreier class inducing this extension.\\

  \noindent As remarked in \refSection{ArtinSchreierTheorySection}, we see that we can even choose $Q$ to be inside the subring
  \begin{align*}
   k[t^{-1}] \subset \localPowerSeriesField{k}{t} = L_0
  \end{align*}
  of polynomials in the variable $t^{-1}$\!. The comments of \refSection{ArtinSchreierTheorySection} also allow us to find a $Q$ with only odd degree terms. We have proved the lemma.
\end{proof}

\noindent The following lemma will be useful later on and we shall leave the proof to the reader.
\begin{lemma}\labelLemma{ClassificationdifferentCor}
 We use the notations of \refThm{ClassificationCor}. If the degree of $Q$ is denoted by $d$ for some odd integer $d$\!, then the degree of local different of 
  \begin{align*}
    L_1 = L^{\left\langle a^2,b\right\rangle} \subset L^{\left\langle b \right\rangle}
  \end{align*}
  is exactly the maximum$\max (4,2d)$\!.
\end{lemma}
\begin{proof}
 One uses the relation (\ref{tvRelationTwo}) together with the Artin--Schreier equation (\ref{ClassificationASEquation}) for the field extension $L_1 \subset L^{\left\langle b \right\rangle}$\!.
\end{proof}

\begin{remark}\labelRemark{distinctionRemark}
 Recall (\refThm{ClassificationCor}) that the extension $L_1 \subset L^{\left\langle b \right\rangle}$ is given by (\ref{ClassificationASEquation}). The proof of \refLemma{ClassificationdifferentCor} shows that if $\degree(Q) \leq 1$\!, then the term $\eta^2 t^{-1}v^{-1}$ of (\ref{ClassificationASEquation}) dominates the degree of different of $L_1 \subset L^{\left\langle b \right\rangle}$\!, i.e. it is then $4$\!. If $\degree(Q) \geq 3$\!, then the term $Q$ dominates this. In \refSection{ProofSection}, we shall prove that all supersimple actions lift to characteristic $0$\!. There we shall distinguish a supersimple action according to the distinction remarked here, i.e. according to the degree of different of  $L_1 \subset L^{\left\langle b \right\rangle}$\!, and we shall need to adapt our lifting technique according to the case we are considering.
\end{remark}
\section{Good reduction of Galois closures}\labelSection{GoodReductionOfGaloisClosures}
\noindent Before we give a brief introduction and overview on this section, we first set some notation. Let \begin{align*} \label{FirstExtension} C_1 \rightarrow \Proj{K} =:C_0,\quad C_2 \rightarrow C_1\end{align*} be $\cyclic{2}$\!-Galois covers of smooth projective $K$\!-curves. We shall assume that the composite extension $ C_2 \rightarrow C_0 \simeq \Proj{K}$ of degree $4$ is not a Galois cover. We let $C_3 \rightarrow C_0 \simeq \Proj{K}$ be the Galois closure.\\

\noindent In Section \ref{GalTheoryDFour} it was shown that  we can identify the Galois group $\Gal{C_3}{\Proj{K}}$ with $\dFour$ in such a manner that $C_2$ is the quotient of $C_3$ under the subgroup \begin{align*}\label{subgroup}\left\langle b \right\rangle \subset \dFour \simeq \Gal{C_3}{C_0}.\end{align*} From now on we shall assume this to be the case. \\

\noindent In this section we shall be concerned with the following question: what reduction conditions on the intermediate cover $C_2 \rightarrow C_0 \simeq \Proj{K}$ are necessary to conclude that the curve $C_3$ has good reduction? In \refSection{MainSection}, we shall make specific choices for the curves $C_0$\!,$C_1$ and $C_2$ which will satsify these conditions. These choices will be such that after studying their reductions, we shall show that by localizing and completing these covers at their branch points, we obtain lifts for all supersimple $\dFour$\!-actions. For our purposes it is convenient to assume that $g(C_2) \geq 1$. In this section we shall place no restrictions on the genus of $C_1$\!, however, in \refSection{MainSection} we shall work only with the case that $g(C_1) = 0$\!, i.e. $C_1 \simeq \Proj{K}$\!.\\

\noindent One sees that if $C_3$ has potentially good reduction, then so must the curve $C_2$. Therefore, we shall always assume that $C_2$ admits a smooth model ${\cal{C}}_2$. We introduce the following assumption on the cover $C_2 \rightarrow C_0$\!.
\begin{assumption}[`Good reduction' Assumption]\labelAssumption{smoothnessAssumption}
  There exists smooth models ${{\cal{C}}}_i$\!, $i = 0,1$\!, of the curves $C_i$\!, $i = 0,1$\!, together with {\emph{finite}} maps \begin{align}
   {\cal{C}}_2 \rightarrow {\cal{C}}_1 \rightarrow {\cal{C}}_0
  \end{align}
  which have generc fibre $C_2 \rightarrow C_1 \rightarrow C_0$. The induced map of smooth $k$\!-curves
  \begin{align}\label{specialFibreCovering}
   {{\cal{C}}}_{2,k} \rightarrow {{\cal{C}}}_{1,k} \rightarrow {{\cal{C}}}_{0,k} \simeq \Proj{k}
  \end{align}
  is a {\emph{separable}} cover of degree $4$\!. Furthermore, we assume that (\ref{specialFibreCovering}) is {\emph{totally}} branched at some point $x \in {\cal{C}}_{0,k}$\!.
\end{assumption}
\begin{remark}
 It follows from Liu--Lorenzini \cite{LiuLorenzini} Proposition 1.6 that since ${{\cal{C}}}_2$ is smooth, the quotients ${{\cal{C}}}_1$ and ${{\cal{C}}}_0$\! are also smooth $R$\!-curves. Furthermore, it follows from Liu \cite{Liu} Proposition 10.3.38 that if furthermore $g(C_1) \geq 2$\!, then (\ref{specialFibreCovering}) is separable.
\end{remark}
\noindent Furthermore, one sees that the following assumption, which does not necessarily hold, is necessary to deduce potentially good reduction for the curve $C_3$\!.
\begin{assumption}[`NonGalois reduction' Assumption]\labelAssumption{nonGaloisAssumption}
  The special fibre cover (\ref{specialFibreCovering}) \begin{align*}   {{\cal{C}}}_{2,k} \rightarrow {{\cal{C}}}_{1,k} \rightarrow {{\cal{C}}}_{0,k} \simeq \Proj{k}
  \end{align*} is not Galois.
\end{assumption}

\noindent Let us now study the stable model $\hat{{{\cal{C}}}}_3$ of the curve $C_3$\!. The group $\dFour$ acts on this model, and we denote by $\hat{{{\cal{C}}}}_i$\!, $i \in \{0,1,2\}$\!, the quotients of this model corresponding to the $K$\!-curves $C_0$\!, $C_1$ and $C_2$ respectively. It is known that all these are themselves semistable $R$\!-curves, see for instance Raynaud \cite{RaynaudFest} Appendice.\\

\noindent Since ${{\cal{C}}}_2$ is a smooth $R$\!-curve with positive genus, we see that there exists a birational blowup morphism
\begin{align}\label{blowupMorphismSecond}
  \hat{{\cal{C}}}_2 \rightarrow {{\cal{C}}}_2.
 \end{align}
Therefore, we may conclude by the universal property of quotient schemes, that similar blowup morphisms
\begin{align}\label{blowupMorphism}
  \hat{{\cal{C}}}_1 \rightarrow {{\cal{C}}}_1,\quad \hat{{\cal{C}}}_0 \rightarrow {{\cal{C}}}_0
 \end{align}
 exist for ${{\cal{C}}}_1$ and ${{\cal{C}}}_0$\!, even if their genera are $0$.\\

\noindent We denote the strict transform of the smooth $k$\!-curve ${{\cal{C}}}_{2,k}$ under the map of (\ref{blowupMorphismSecond}) by $\Gamma_2$\!, and using (\ref{blowupMorphism}) we define the components $\Gamma_1$ and $\Gamma_0$ similarly. Each $\Gamma_i$\!, for $i = 0,1,2$\!, is therefore a smooth $k$\!-curve, and furthermore, we have a separable degree-\!$4$ covering 
\begin{align}
  \Gamma_2 \rightarrow \Gamma_1 \rightarrow \Gamma_0
\end{align}
which is nothing else than the covering
\begin{align}
  {\cal{C}}_{2,k} \rightarrow {\cal{C}}_{1,k} \rightarrow {\cal{C}}_{0,k} \simeq \Proj{k}.
\end{align}
\noindent Let $\Gamma_3$ be any component of $\hat{{\cal{C}}}_{3,k}$ which maps surjectively onto $\Gamma_2$ under the finite map $\hat{{{\cal{C}}}}_{3,k} \rightarrow \hat{{{\cal{C}}}}_{2,k}$\!.  Since $\hat{{{\cal{C}}}}_{3,k}$ was assumed to be the stable model of $C_3$\!, we see that each component of $\hat{\cal{C}}_{3,k}$ is reduced, and, in particular, the closed subscheme $\Gamma_3$ is an integral scheme. We may therefore consider the extension of function fields
\begin{align}\label{FieldExtensionTower}
  k(\Gamma_0) \subset k(\Gamma_1) \subset k(\Gamma_2) \subset k(\Gamma_3).
\end{align}
\begin{prop}\labelProp{DecompProp}
  The component $\Gamma_3$ is the only component of $\hat{\cal{C}}_{3,k}$ mapping surjectively onto $\Gamma_2$\!. Furthermore, the field extension (\ref{FieldExtensionTower}) is a $\dFour$\!-Galois extension. 
\end{prop}
\begin{proof}
 Let $D(\Gamma_3)$ (respectively $I(\Gamma_3)$\!) denote the decomposition (respectively inertia) group of $\Gamma_3$\!. Let $L$ be the separable closure of $k(\Gamma_0)$ inside the normal field extension (\ref{FieldExtensionTower}). There exists an exact sequence of groups (see Serre \cite{SerreLocalFields} Proposition I.20)
 \begin{align*}
  0 \rightarrow I(\Gamma_3) \rightarrow D(\Gamma_3) \rightarrow \Gal{L}{k(\Gamma_0)} \rightarrow 0.
 \end{align*}
 Notice that by \refAssumption{smoothnessAssumption} the Galois extension $k(\Gamma_0) \subset L$ contains the subextension $$k(\Gamma_0) \subset k(\Gamma_2).$$ Furthermore, by \refAssumption{nonGaloisAssumption}, we see that $[L : k(\Gamma_0)] > 4$ and therefore, the order of $\Gal{L}{k(\Gamma_0)}$ must exceed $4$\!. However, $D(\Gamma_3)$ is a subgroup of $\dFour$\!, and therefore, the result follows.
\end{proof}

\noindent Our next step is to study the normalization of the component $\Gamma_3$\!. Let $\tilde{\Gamma}_3$ denote the normalization of $\Gamma_3$\!. In order to deduce smoothness of the stable $R$\!-curve $\hat{\cal{C}}_3$\!, we shall now ask for a condition under which the geometric genus $g(\tilde{\Gamma}_3)$ of $\tilde{\Gamma}_3$ is equal to the geometric genus $g(C_3)$ of the generic fibre $C_3$\!. It is known that the latter is never strictly less than the former. Furthermore, since $\hat{\cal{C}}_3$ is assumed to be the stable model of the $K$\!-curve $C_3$\!, equality of $g(\tilde{\Gamma}_3)$ and $g(C_3)$ would imply smoothness of $\hat{\cal{C}}_3$\!. We thus proceed to bounding $g(\tilde{\Gamma}_3)$ from below.

\begin{assumption}[`Different' Assumption]\labelAssumption{differentAssumption}
  We assume the degree of geometric different of the cover $ C_3 \rightarrow C_2 $ is $2$\!.
\end{assumption}
\begin{lemma}
  The genus of $C_3$ is $2g(C_2)$\!. In particular, we have the following inequalities.
  \begin{align}
    g(\tilde{\Gamma}_3) \leq 2g(C_2).
  \end{align}
\end{lemma}
\begin{proof}
 Apply the Hurwitz Formula to the cover of smooth $K$\!-curves $C_3 \rightarrow C_2$\!.
\end{proof}
\begin{thm}\labelThm{SmoothnessThm}
  The curve $\hat{{\cal{C}}}_3$ is a smooth $R$\!-curve.
\end{thm}
\begin{proof}
 \noindent By \refAssumption{differentAssumption}, the cover of $K$\!-curves $C_3 \rightarrow C_2$ has exactly two geometric branch points $x_1, x_2$\!, and after possibly extending $K$\!, we may assume that these two points are distinct points of $C_2(K)$\!. From Theorem 1 of Sa\"{i}di \cite{Saidi}, we see that since $\Gamma_3 \rightarrow \Gamma_2$ is a separable covering, both $x_1$ and $x_2$ specialize to the same point $x$ of ${\cal{C}}_{2,k}$\!. Note that $\tilde{\Gamma}_3 \rightarrow \Gamma_2$ is branched at this point. This implies that \begin{align}g(\tilde{\Gamma}_3) \geq 2g(\Gamma_2) = 2g(C_2),\end{align} and hence is equal to exactly this. Therefore, $\hat{\cal{C}}_3$ is a smooth $R$\!-curve.
\end{proof}
\section{Lifting supersimple $\dFour$\!-actions.}\labelSection{MainSection}
In this section we shall give a method for producing covers of curves which satisfy the assumptions  needed to apply the results in the previous section. Let us first set some notation and then we explain our goals and strategy.\\

\noindent In \refSection{GoodReductionOfGaloisClosures}, we dealt with towers of $\cyclic{2}$\!-covers $C_2 \rightarrow C_1 \rightarrow C_0$ of composite degree $4$\!. Our first step is to construct suitable choices for the curves $C_1$ and $C_0$\!. We let ${{\cal{C}}}_0$ denote the projective $R$\!-line $\Proj{R,t}$ with parameter $t$\!. \\

\noindent To define the $R$\!-curve ${\cal{C}}_1$\!, we define an algebraic extension of $K(t)$ by adjoining the element $v$\!, where $v$ satisfies the relation
 \begin{align*}
   t^{-1} = v^{-2} - v^{-1}.
 \end{align*}
\noindent We now define ${\cal{C}}_1$ to be the normalization of the projective line ${\cal{C}}_0$ inside the field $K(t)(v)$\!. We leave for the reader to verify that ${\cal{C}}_1$ is again a projective $R$\!-line with parameter $v$\!, and that the induced special fibre cover ${{\cal{C}}}_{1,k} \rightarrow {{\cal{C}}}_{0,k}$ is a separable cover of smooth projective lines. By localizing and completing at the point $t = 0$\!, we see that the $\cyclic{2}$\!-Galois cover \begin{align*}\Proj{R,v} = {\cal{C}}_1 \rightarrow {\cal{C}}_0 = \Proj{R,t}\end{align*} already provides a lift for the $\cyclic{2}$\!-Galois extension of local fields $L_0 \subset L_1$ of Notation \ref{FieldNotation}. \\

\noindent Now we want to construct some $\cyclic{2}$\!-Galois extensions of the curve ${\cal{C}}_1$\!.  Let $F$ and $G$ be two elements of $R[v^{-1}] \subset K(v^{-1})$\!. We denote the reductions of $F$ and $G$ to the ring $k[v^{-1}]$ by $\overline{F}$ and $\overline{G}$\!, respectively. We define a field extension $K(v) \subset K(v,w)$ where $w$ satisfies
\begin{align}\label{defEq}w^2 - wG = F.\end{align} We let ${\cal{C}}_2^{F,G}$ be the normalization of ${\cal{C}}_1 = \Proj{R,v}$ inside $K(v,w)$\!. We have included the superscripts $F$ and $G$ to emphasize that our definition depends on the choices of $F$ and $G$\!. \\

\noindent Our strategy now is to find suitable $F$ and $G$ such that the generic fibre of the finite tower of $\cyclic{2}$\!-Galois extensions
 \begin{align}
   {\cal{C}}_2^{F,G} \rightarrow {\cal{C}}_1 = \Proj{R,v} \rightarrow {\cal{C}}_0 = \Proj{R,t}
 \end{align}
satisfies Assumptions \ref{thm : smoothnessAssumption}, \ref{thm : nonGaloisAssumption} and \ref{thm : differentAssumption}. \\

\noindent To check the `good reduction' (\refAssumption{smoothnessAssumption}), the form of equation (\ref{defEq}) will be useful (\refLemma{ConstructionLemma}). However, to check \refAssumption{differentAssumption}, we shall need to rewrite this equation in a Kummer form. Here we shall restrict the choices of $F$ and $G$ (\refLemma{GoodFormLemma}). A further restriction (Lemmas \ref{thm : HelpfulGoodFormLemmaSmall} and \ref{thm : HelpfulGoodFormLemma}) on the choices of $F$ and $G$ will also aid us in checking that the `reduction is not Galois' (\refAssumption{nonGaloisAssumption}).\\

\noindent Assume that the degrees of $F$ and $\overline{F}$ are both $2g + 1$\!, where $g$ is some positive integer. Furthermore, assume that the degree of $G$ does not exceed $2g$\!, and that the reduction $\overline{G}$ is a unit of $k$ (and hence of degree $0$\!, but that $\overline{G} \neq 0$ inside $k[v^{-1}]$\!). \\
\begin{lemma}\labelLemma{ConstructionLemma}

{\noindent \textnormal(a)} The scheme ${\cal{C}}_2^{F,G}$ is a smooth projective $R$\!-curve of genus $g$\!. Furthermore, the action of the Galois group $\Gal{K(v,w)}{K(v)}$ extends to the scheme ${\cal{C}}_2^{F,G}$\!, and the quotient of ${\cal{C}}_2^{F,G}$ by this action is ${\cal{C}}_1 = \Proj{R,v}$\!. Lastly, the induced map of special fibres
 \begin{align*}
   {\cal{C}}_{2,k}^{F,G} \rightarrow {\cal{C}}_{1,k} \simeq \Proj{k,v}
 \end{align*}
 is generically separable, and is in fact branched uniquely at the point $v = 0$\!. 

{\noindent \textnormal(b)} By localizing and completing at this point, the cover ${{\cal{C}}}_2^{F,G} \rightarrow {\cal{C}}_1$ induces a $\cyclic{2}$\!-Galois extension of $\localPowerSeriesField{k}{v}$ generated by $w$\!, where $w$ satisfies
 \begin{align}
  w^2 - w\overline{G} = \overline{F}.
 \end{align}
 The local degree of different is $2g+2$\!.
\end{lemma}
\begin{proof}
 This is essentially Exercise 10.1.9 of Liu \cite{Liu}.
\end{proof}

\noindent  So far, we have constructed a tower of smooth projective curves
\begin{align}\label{RCovering}
  {\cal{C}}_{2}^{F,G} \rightarrow {\cal{C}}_1 \rightarrow {\cal{C}}_0 \simeq \Proj{R}.
\end{align} 
For convenience we set $C_2^{F,G}:={\cal{C}}_{2,K}^{F,G}$ and similarly for $C_1$ and $C_0$\!. As in \refSection{GoodReductionOfGaloisClosures}, we define $C_3^{F,G}$ to be the Galois closure of $C_2^{F,G} \rightarrow C_0$\!. Notice that in order to apply the results of \refSection{GoodReductionOfGaloisClosures}, we also need to to know that $C_2^{F,G} \neq C_3^{F,G}$\!. This will be true for the choices of $F$ and $G$ that we shall later choose. 
\begin{lemma}\labelLemma{GoodFormLemma}
   Assume that we can find an element \begin{align*}H \in R[t^{-1}] \subset K(t) \subset K(v),\end{align*} as well as an $\eta \in R$\! such that the following identity holds.
  \begin{align*}
    4F + G^2 = (1 - 2\eta v^{-1})H.
  \end{align*}
 If $\eta \neq 1$\!, then $C_2^{F,G} \neq C_3^{F,G}$\!, i.e. the cover $C_2^{F,G} \rightarrow C_0$ is not Galois, and furthermore, the degree of geometric different of $C_3^{F,G} \rightarrow C_2^{F,G}$ is $2$\!.
\end{lemma}
\begin{proof}
 By construction, the cover \begin{align*}\Proj{K,v} \simeq C_1 \rightarrow C_0\simeq \Proj{K,t}\end{align*} is ramified at exactly $v = 0$ and $v = 2$\!. Notice that the function field $K(v,w)$ of $C_2^{F,G}$ is also generated over $K(v)$ by $w'$\!, where $w' = -2w + G$\!, i.e. $w'$ satisfies the following Kummer equation
\begin{align*}
 (w')^2 = (-2w + G)^2 = 4F + G^2 = (1 - 2\eta v^{-1})H.
\end{align*}

\noindent Therefore, the cover $C_2^{F,G} \rightarrow C_1$ is branched at exactly $v = 0$\!, $v = 2\eta$ and the zeros  of $$H \in K(t) \subset K(v).$$ If $\eta \neq 1$\!, then the conjugate of the point $v = 2\eta$ (under the action of $\Gal{C_1}{C_0}$\!)  is not branched in the cover $C_2^{F,G} \rightarrow C_1$\!. This already implies that $C_2^{F,G} \rightarrow C_0$ is not Galois, i.e. $C_2^{F,G} \neq C_3^{F,G}$\!. Furthermore, one checks that the points of $C_2^{F,G}$ lying above the conjugate of the point $v = 2\eta$ are exactly the branch points of the cover $C_3^{F,G} \rightarrow C_2^{F,G}$\!. There are exactly two points, and hence the degree of geometric different of $C_3^{F,G} \rightarrow C_2^{F,G}$ is $2$\!.
\end{proof}
Before we state the main theorems of this section, we give two computational results. Let $L_0 \subset L$ be a local supersimple $\dFour$\!-Galois extension in characteristic $2$\!. Recall from \refLemma{ClassificationdifferentCor} and \refRemark{distinctionRemark} that we can distinguish between two cases, namely the case where the local different degree of $L^{\left\langle a^2,b \right\rangle} \subset L^{b}$ is $4$\!, and the case where it is $2d$\!, for some odd integer $d$\!. In proving that all supersimple actions lift to characteristic $0$\!, we shall deal with these two cases seperately. In both cases, we shall need a similar computation, and it is these that we state in the following two lemmas. Both of these results are essentially computations, and we used the computer package Magma to verify our calculations.
\begin{lemma}\labelLemma{HelpfulGoodFormLemmaSmall}
 Let $\eta \in R^{*}$\!. We assume $R$ has been extended to include a solution, $\beta$\!, of the following equation.
 \begin{align}\label{HGFLE2}
   \beta^2 + \sqrt{2}\beta + \eta = 0.
 \end{align}
 Let $Q' \in R[t^{-1}]$ of degree less than or equal to $1$\!. Then we have the following identity.
 \begin{align*}
  (1 - 2\eta v^{-1})(1 + 2\beta^2t^{-1} + 4Q') = G^2 + 4F,
 \end{align*}
 where 
 \begin{align*}
  G:=1 + \sqrt{2}\beta  v^{-1} 
 \end{align*}
 and 
 \begin{align*}
  F:=Q' - \eta \beta^2 v^{-1} t^{-1} - 2\eta v^{-1} Q'.
 \end{align*}
\end{lemma}
\begin{lemma}\labelLemma{HelpfulGoodFormLemma}
 Let $\eta$ and $\beta$ be as in \refLemma{HelpfulGoodFormLemmaSmall}. Let $m$ be a positive integer, and let $Q'$ be any element of $R[t^{-1}]$ of degree {\emph{strictly}} less than $2m$\!. Furthermore, let $\gamma \in R^{*}$ be any unit of $R$\!. Then we have the following identity.
 \begin{align*}
  (1 - 2\eta v^{-1})(1 + 2\beta^2t^{-1} + 2\gamma^2t^{-2m} + 2\sqrt{2}\gamma t^{-m} + 4Q') = G^2 + 4F,
 \end{align*}
 where 
 \begin{align*}
  G:=1 + \sqrt{2}\beta  v^{-1} + \sqrt{2}\gamma t^{-m}
 \end{align*}
 and 
 \begin{align*}
  F:=Q' - \eta \beta^2 v^{-1} t^{-1} - 2\eta v^{-1} Q' - \eta \gamma^2 t^{-2m} v^{-1} - \sqrt{2}\eta\gamma v^{-1} t^{-m} - \gamma \beta v^{-1} t^{-m}.
 \end{align*}
\end{lemma}
\begin{remark}\labelRemark{HGFLERemark}
 The equation (\ref{HGFLE2}) implies that we have the following equality in $k$ after reduction
 \begin{align} \overline{\beta}^2 = \overline{\eta}.\end{align}
\end{remark}
\begin{remark}\labelRemark{reductionFRemark}
  The polynomial $F$ in \refLemma{HelpfulGoodFormLemma} reduces to 
  \begin{align*}
    \overline{F}:=\overline{Q'} + \eta\beta^2 v^{-1} t^{-1}. 
  \end{align*}
  If $F$ is selected as in \refLemma{HelpfulGoodFormLemma}, then  it reduces to
  \begin{align*}
    \overline{F}:=\overline{Q'} + \eta\beta^2 v^{-1} t^{-1} + \eta\gamma^2 t^{-2m}v^{-1} + \beta\gamma t^{-m} v^{-1}.
  \end{align*}
  In both cases $G$ reduces to the constant polynomial $1 \in k$\!.
\end{remark}
\noindent The following theorem is our first main result. It constructs a family of $\dFour$\!-Galois covers which, by localizing and completing at branch points, induce local supersimple extensions after reduction. We define the normal $R$\!-curve ${\cal{C}}_3^{F,G}$ to be the normalization of ${\cal{C}}_2^{F,G}$ inside the extension $C_3^{F,G} \rightarrow C_2^{F,G}$.
\begin{thm}\labelThm{FGgoodreduction}
 Let either $\eta$\!, $Q'$ be as in \refLemma{HelpfulGoodFormLemmaSmall}, or let $\eta$\!, $\gamma$\!, $m$ and $Q'$ be as in \refLemma{HelpfulGoodFormLemma}, and let $F$ and $G$ be selected as in these lemmas. Consider the $\dFour$\!-Galois extension of normal projective $R$\!-schemes
 \begin{align}\label{consideredExtension}
   {\cal{C}}_3^{F,G} \rightarrow {\cal{C}}_2^{F,G} \rightarrow {\cal{C}}_1 \rightarrow {\cal{C}}_0 = \Proj{R,t}.
 \end{align}
 Then each ${\cal{C}}_{i}^{F,G}$ is a smooth $R$\!-scheme. Furthermore, by localizing and completing at the point $t = 0$ of the scheme ${\cal{C}}_0 = \Proj{R,t}$\!, we obtain a lifting of the local $\dFour$\!-Galois extension obtained by taking the Galois closure of 
 \begin{align*}
  \localPowerSeriesField{k}{t} \subset \localPowerSeriesField{k}{v} \subset \localPowerSeriesField{k}{v}(w),
 \end{align*}
 where $w$ satisfies
 \begin{align*}
   w^2 - w = \overline{F}. 
 \end{align*}
 Here $\overline{F}$ denotes the reduction of the polynomial $F$\!, refer to \refRemark{reductionFRemark} for an explicit expression of $\overline{F}$.
\end{thm}
\begin{proof}
 We shall proof the theorem in the case that $F$ and $G$ have been selected as in \refLemma{HelpfulGoodFormLemma} and leave the (easier) case of \refLemma{HelpfulGoodFormLemmaSmall} to the reader.\\

 \noindent First we see from \refLemma{ConstructionLemma} that \refAssumption{smoothnessAssumption} is satisfied for the extension (\ref{consideredExtension}). In fact, the model ${\cal{C}}_2^{F,G}$ is a smooth model for its generic fibre $C_2^{F,G}$\!, and by construction the special fibre subcover \begin{align}{\cal{C}}_{2,k}^{F,G} \rightarrow {\cal{C}}_{1,k} = \Proj{k,v} \rightarrow {\cal{C}}_0 = \Proj{k,t}\end{align} is separable. \refLemma{GoodFormLemma} tells us that \refAssumption{differentAssumption} is also satisfied for this extension. \\

 \noindent Let us check that the induced cover \begin{align}\label{compositeGaloisCover}{{\cal{C}}}_{2,k}^{F,G} \rightarrow {{\cal{C}}}_{1,k} \rightarrow {{\cal{C}}}_{0,k}\end{align} is not a Galois cover, thereby verifying \refAssumption{nonGaloisAssumption}. By localizing and completing at the point $v = 0$ of ${{\cal{C}}}_{1,k} \simeq \Proj{k,v}$\!, we obtain a cover of $\localPowerSeriesField{k}{v}$ generated by $w$\!, where $w$ satisfies
 \begin{align*}
   w^2 - w =  \overline{F} = {\eta}^2t^{-1}v^{-1} + Q' + \eta\gamma^2t^{-2m}v^{-1} + \beta\gamma t^{-m}v^{-1}.
 \end{align*}
 \noindent One checks that the composite field extension $\localPowerSeriesField{k}{t} \subset \localPowerSeriesField{k}{v} \subset \localPowerSeriesField{k}{v}(w)$ is not Galois if $\overline{\eta} \neq 1$\!. Therefore, the composite cover (\ref{compositeGaloisCover}) cannot be Galois. By \refThm{SmoothnessThm}, we see that the curve $C_3^{F,G}$ has potentially good reduction. Since the smooth model ${\cal{C}}_2^{F,G}$ of $C_2^{F,G}$ is unique (recall that $g(C_2) \geq 1$), we see that ${\cal{C}}_3^{F,G}$ is smooth.  We are done.
\end{proof}

\section{Proof of main result}\labelSection{ProofSection}
\noindent The aim of this section is to prove our main result.
\begin{thm}\labelThm{MainThm}
 All supersimple $\dFour$\!-actions lift to characteristic $0$\!.
\end{thm}
 \noindent Assume throughout this section that we have been given a supersimple $\dFour$\!-Galois extension of local power series fields  \begin{align*}\label{givenSupersimpleExtension} L_0:=\localPowerSeriesField{k}{t} \subset L.\end{align*}  We use the notation of Sections \ref{thm : SupersimpleSection} and \ref{thm : ClassificationSection}. In particular, we set $L_1:=L^{\left\langle a^2,b \right\rangle}$ with parameter $v$\!, where $v$ and $t$ are related by 
\begin{align}
 v^{-2} - v^{-1} = t^{-1}.
\end{align}

 \noindent We have already pointed out (\refRemark{UniqueRemark}) that the field extension $L_0 \subset L$ is completely determined by the subextension \begin{align*}L_0 \subset L_1 \subset L_2:=L^{\left\langle b \right\rangle}.\end{align*} By \refLemma{ClassificationdifferentCor}, we see that there are two cases to consider, namely the case that the degree of different of $L_1 \subset L_2$ is $4$ or the case that it is $2d$\!, where $d > 1$ is an odd integer. \\

\noindent In both cases, we apply \refThm{FGgoodreduction} for suitable choices of $\eta$\!, $\gamma$\!, $m$ and $Q'$\!. In the first case, we shall choose the $F$ and $G$ as in \refLemma{HelpfulGoodFormLemmaSmall}, and in the second as in \refLemma{HelpfulGoodFormLemma}. We shall give the details only for the second case, and leave the detailed proof of the first (easier) case to the reader.

\begin{proofLater}
 \noindent We assume that the local different degree of $L_1 \subset L_2$ is of the form $2d$\!, where $d > 1$ is an odd integer. Define $m$ by the relation $d = 2m+1$\!.\\

 \noindent From \refThm{ClassificationCor} and \refLemma{ClassificationdifferentCor}, we see that we can find a polynomial $Q \in k[t^{-1}]$ of degree exactly $2m+1$\!, as well as an $\eta \in k-\finiteField{2}$\!, such that the extension $L_1 \subset L_2$ is generated by $w$\!, where $w$ satisfies 
\begin{align}\label{LEXT}
 w^2 - w = \eta^2 t^{-1}v^{-1} + Q.
\end{align}

 \noindent Recall (\refThm{ClassificationCor}) that we can choose $Q$ to have only odd degree terms in $t^{-1}$\!. Furthermore, since $k$ was assumed algebraically closed, we can find a $\gamma \in k^{*}$ such that 
 \begin{align*}
   Q = \gamma^2 \eta t^{-(2m+1)} + Q',
 \end{align*}
 where $Q' \in k[t^{-1}]$ has degree strictly smaller than $2m$\!. \\

 \noindent Let us lift the elements $\eta$ and $\gamma$ to units of $R$\!. We abuse notation and denote these lifts again by $\eta$ and $\gamma$\!, respectively. We then choose a polynomial $Q' \in R[t^{-1}]$\!, of degree less than $2m$\!, which lifts the polynomial $Q' \in k[t^{-1}]$\!. \\

 \noindent We now apply \refThm{FGgoodreduction} with these choices of $\eta$\!, $\gamma$\!, $m$ and $Q'$ and we choose $F$ and $G$ as in \refLemma{HelpfulGoodFormLemma}. In view \refThm{FGgoodreduction}, we only need to check  that the Artin--Schreier class of $\overline{F}$ is the same as $\eta^2 t^{-1}v^{-1} + Q$\!. \\

 \noindent From \refRemark{reductionFRemark}, we have the following equality of Artin--Schreier classes inside $\cohH^1(G_{L_1},\cyclic{2})$
 \begin{align*}
   \class{\overline{F}} &= \class{\eta\beta^2 v^{-1} t^{-1} + {Q'} + \eta\gamma^2 t^{-2m}v^{-1} + \beta\gamma t^{-m} v^{-1}}.
 \end{align*}
 Since $\beta^2 = \eta$ inside $k$\! (\refRemark{HGFLERemark}), we see that 
 \begin{align*}
  \class{\overline{F}} = \class{\eta^2 v^{-1} t^{-1} + {Q'} + \eta\gamma^2t^{-2m - 1}} 
			= \class{\eta^2v^{-1}t^{-1} + Q}.
 \end{align*}
 This is exactly the class of the extension $L_1 \subset L_2$\!, see (\ref{LEXT}). We conclude by applying \refThm{FGgoodreduction}.
\end{proofLater}

\end{document}